\documentclass[11pt,twoside,letter]{article}
\usepackage{amsthm, amsmath, amssymb, amsfonts, graphicx, epsfig}
\usepackage{accents}

\usepackage{bbm}
\usepackage{mathrsfs}

\usepackage{epstopdf}

\usepackage{graphicx}
\usepackage[font=sl,labelfont=bf]{caption}
\usepackage{subcaption}

\usepackage{enumerate}

\usepackage[colorlinks=true, pdfstartview=FitV, linkcolor=blue,
citecolor=blue, urlcolor=blue]{hyperref}
\usepackage[usenames]{color}

\usepackage{url}
\makeatletter\def\url@leostyle{%
	\@ifundefined{selectfont}{\def\UrlFont{\sf}}{\def\UrlFont{\scriptsize\ttfamily}}} \makeatother

\usepackage[margin=1.0in, letterpaper]{geometry}

\newtheorem{theorem}{Theorem}
\newtheorem{proposition}[theorem]{Proposition}
\newtheorem{lemma}[theorem]{Lemma}

\theoremstyle{definition}

\newtheorem{example}[theorem]{Example}
\theoremstyle{remark}
\newtheorem{remark}[theorem]{Remark}

\numberwithin{equation}{section}
\numberwithin{theorem}{section}

\def\cB{\mathcal{B}}
\def\cC{\mathcal{C}}
\def\cD{\mathcal{D}}

\def\cN{\mathcal{N}}

\def\bE{\mathbb{E}}
\def\bF{\mathbb{F}}

\def\bN{\mathbb{N}}

\def\bP{\mathbb{P}}

\def\bR{\mathbb{R}}

\def\sF{\mathscr{F}}

\newcommand{\set}[1]{\{#1\}}            
\newcommand{\norm}[1]{ \| #1 \| }       

\DeclareMathOperator{\dif}{d \!}        
\DeclareMathOperator{\wlim}{{w}-\lim}        

\title{ \vspace{-3em} 
Statistical analysis of discretely sampled semilinear SPDEs: a power variation approach}

\makeatletter
\def\and{%
\end{tabular}%
\begin{tabular}[t]{c}}%
\def\@fnsymbol#1{\ensuremath{\ifcase#1\or a\or b\or c\or
		d\or e\or f\or g\or h\or i\else\@ctrerr\fi}}
\makeatother

\author{
	Igor Cialenco\,\thanks{Department of Applied Mathematics, Illinois Institute of Technology
		\newline \hspace*{1.45em}  10 W 32nd Str, Building RE, Room 220, Chicago, IL 60616, USA
		\newline \hspace*{1.45em}  Emails: \url{cialenco@iit.edu}, URL: \url{http://cialenco.com}
		\vspace{0.5em}} ,
	\and
	Hyun-Jung Kim\,\thanks{Department of Mathematics, University of California Santa Barbara, 
		\newline \hspace*{1.45em} South Hall, Room 6607,
		University of California Santa Barbara, CA 93106, USA
		\newline \hspace*{1.45em} Email: \url{hjkim@ucsb.edu},
 URL: \url{https://sites.google.com/view/hyun-jungkim}
		\vspace{0.5em}},
	\and
	Gregor Pasemann\,\thanks{Institut f\"ur Mathematik, Technische Universit\"at Berlin, Stra{\ss}e des 17. Juni 136, Berlin, Germany
		\newline \hspace*{1.45em} Email: \url{pasemann@math.tu-berlin.de} 
		\vspace{0.5em}}
}

\date{ {\small 
		First Circulated and this version: March 4, 2021\\
}}

\newcommand{\N}{\mathbb{N}}
\newcommand{\Z}{\mathbb{Z}}
\newcommand{\R}{\mathbb{R}}
\newcommand{\unitvar}{\nu}

\begin{document}

	\maketitle

	\vspace{-2em}
	
	
	{\footnotesize
		\begin{tabular}{l@{} p{350pt}}
			\hline \\[-.2em]
			\textsc{Abstract}: \ & 
			
Motivated by problems from statistical analysis for discretely sampled SPDEs, first we derive central limit theorems for higher order finite differences applied to stochastic process with arbitrary finitely regular paths.  These results are proved by using the notion of $\Delta$-power variations, introduced herein, along with the H\"older-Zygmund norms. Consequently, we prove a new central limit theorem for   $\Delta$-power variations of the iterated integrals of a fractional Brownian motion (fBm). These abstract results, besides being of independent interest, in the second part of the paper are applied to estimation of the drift and volatility coefficients of semilinear stochastic partial differential equations in dimension one, driven by an additive Gaussian noise white in time and possibly colored in space. In particular, we solve the earlier conjecture from \cite{CialencoKimLototsky2019} about existence of a nontrivial bias in the estimators derived by naive approximations of derivatives by finite differences. We give an explicit formula for the bias and derive the convergence rates of the corresponding estimators. Theoretical results are illustrated by numerical examples.  

\\[0.5em]
\textsc{Keywords:} \ &  statistical inference for SPDEs; CLT for iterative integrals; $\Delta$-power variations; fractional Brownian motion; discrete sampling; semilinear SPDEs\\
\textsc{MSC2010:} \ &  Primary 60F05; Secondary 60H15, 62M05, 62G05 62F12. \\[1em]
\hline
\end{tabular}
	}

\section{Introduction}

The main motivation of this work comes from some open problems in statistical analysis of \textit{discretely sampled} stochastic partial differential equations (SPDEs) of the form 
\begin{align}\label{eq:intro-main-bdd}
\mathrm{d}X_t(x) &= -\theta(-\Delta)^{\alpha/2} X_t(x)\mathrm{d}t + F(X_t(x)) \dif t + \sigma (-\Delta)^{-\gamma}\mathrm{d}W_t(x), \quad t>0,\ x\in [0,1], 
\end{align}
where $\alpha>0,\gamma\geq0$ are given, $\theta, \sigma>0$ are the parameters of interest (unknown to the observer), $W$ is a cylindrical Wiener process on $L^2([0,1])$, and $F$ is a (nonlinear) operator acting on some appropriate Hilbert space.  Most of the existing literature on statistical inference for SPDEs is dedicated to linear SPDEs, i.e. $F=0$, with few exceptions \cite{IgorNathanAditiveNS2010,PasemannStannat2019,ACP2020,PasemannEtAl2020,AltmeyerEtAl2020}. Moreover, majority of works were dedicated to continuous time sampling setup; cf. the survey paper \cite{Cialenco2018}. The parameter estimation problem for  (linear) SPDEs when the solution is discretely sampled in space and/or time component was  addressed systematically only recently by quite different methods, and we refer to \cite{CialencoHuang2017,BibingerTrabs2017,BibingerTrabs2019,Chong2019,CialencoDelgado-VencesKim2019,Chong2019a,KainoMasayukiUchida2019,KhalilTudor2019,KhalilTudor2019b,HildebrandtTrabs2019,CialencoKim2020,ShevchenkoEtAl2020}, and to \cite{PiterbargRozovskii1997,PospivsilTribe2007} for earlier studies. The central theme in these works evolves, in one form or another, around power variations of some relevant stochastic processes, which in turn is strongly related to the regularity properties of the solution. For example, when $\alpha=2$, $\gamma=0$, and $F=0$, one can show that for a fixed $x\in(0,1)$, the paths of the process $X_t(x)$ have continuous versions with H\"older order of continuity $1/4-\varepsilon$, for any $\varepsilon>0$. Consequently, as proved in \cite{CialencoHuang2017}, the fourth power variation is finite and yields consistent and asymptotically normal estimators for $\theta$ and $\sigma$. Similar arguments hold true for solutions of SPDEs when the H\"older order of continuity in space or time component is smaller than one. However, this approach cannot be applied directly to SPDEs with regular paths, and the main goal of this work is to develop new methodologies that can treat such cases. Of course, one should not expect that the solution $X_t(x)$ as function of $t$ will get smoother than the paths of a Brownian motion, i.e. almost 1/2 H\"older continuous. On the other hand, it is known, for example when $F=0$, that for any fixed $t>0$, the solution process $X_t(x), x\in(0,1)$, has almost H\"older $2\gamma+\alpha/2-1/2$ regularity in spatial variable $x$, namely the solution gets smoother the more colored (correlated) in space is the driving noise.  One approach is to take the maximal number of (classical) derivatives in $x$, say $m:=\lfloor 2\gamma+\alpha/2-1/2\rfloor$, and expect that $\partial^m_xX_t(x)$ is equal to a fractional Brownian motion with Hurst parameter $2\gamma+\alpha/2-1/2- m$  plus a smooth process, and apply or adapt the existing results on power variations, for example, from \cite{CialencoHuang2017,KhalilTudor2019,KhalilTudor2019b}. However, from statistical point of view, this assumes that the  process $\partial^m_xX_t(x), \ x\in(0,1)$ is observed, which practically speaking is an unrealistic assumption. One way to overcome this drawback, is to approximate the derivatives by using the discrete measurements of the solution itself, for example by finite differences. However, such approximations typically will yield a nontrivial and non-vanishing bias in the estimators - a phenomena noticed in \cite{CialencoKimLototsky2019} through numerical experiments for SPDEs driven by space-only noise and with $m=1$, and later in \cite{CialencoKim2020} the bias was explicitly given and the asymptotic properties of the estimator were formally proved. We built on these line of ideas, and we focus our study on discretely  sampled (in space) of semilinear SPDEs.  

A key concept of this paper is to track and use the classical regularity of a continuous function in terms of conveniently chosen integro-difference operators, for which we use the H\"older-Zygmund norms and spaces rather than classical H\"older or Sobolev norms and spaces.  To deal with the higher order finite differences and their power variations, we introduce the notion of $\Delta$-power variation,  and prove that the central limit theorems for $\Delta$-power variations are invariant under smooth perturbations; see Section~\ref{sec:PowerVar}. We note that the idea of using quadratic variation of higher order finite differences have been used, for example, in estimation of self-similarity order of self-similar processes; see, for instance, \cite[Section~5.6]{Tudor2013} and references therein. We derive a new central limit theorem for $\Delta$-power variations of iterated integrals of a fractional Brownian motion (fBm) (see Section~\ref{sec:fBM}), where we also explicitly compute the asymptotic variance. These novel results are of independent interest, contributing to the literature on limit theorems for fractional type processes, but in addition, these results provide a method for building consistent and asymptotically normal estimators for discretely sampled process with smooth paths, such as the SPDEs mentioned earlier.

Statistical analysis of semilinear SPDEs is investigated in Section~\ref{sec:SPDE-BDom} and Section~\ref{sec:unboundedDom}.  We study the estimation of the drift $\theta$ and volatility $\sigma$ of \eqref{eq:intro-main-bdd}, under fairly general assumptions on the nonlinear part, assuming that the solution is sampled discretely in the spatial component $x$ at one fixed time instance $t>0$. 
Similarly to the above cited works on nonlinear SPDEs, we first use the so-called splitting of the solution argument, where the solution is written as $X=\overline{X} + \widetilde{X}$, where $\overline X$ is the solution of the linear SPDE and $\widetilde{X}$ solves the corresponding nonlinear random PDE (see equations \eqref{eq:DirichletConditionLinearEquation} and \eqref{eq:barXBoundDom}). Usually $\widetilde X$ is smoother than $\overline{X}$, which allows to argue that the estimation problem can be reduced to the linear case. The latter is reduced to the results on fBm by proving that the highest order (classical) derivative of $\overline{X}$ has the same probability law as a smoothly perturbed fBm. Assuming that one of the coefficients $\sigma$ or $\theta$ is known we derive an estimator for the second coefficient, prove its consistency and provide its rate of convergence.  We note that, the results in \cite{CialencoHuang2017}, which is the closest  in spirit to this manuscript, considers only linear equations driven by space-time white noise, i.e. $\alpha=2$, $\gamma=0$, and $F=0$. The results presented in this manuscript are the first ones on parameter estimation for SPDEs with arbitrarily regular paths that are discretely sampled in physical spatial domain. As a second application of general results of Section~\ref{sec:fBM}, in Section~\ref{sec:unboundedDom} we study parameter estimation problem for a version of SPDEs \eqref{eq:intro-main-bdd} on the whole space. Namely, same as in  \cite{KhalilTudor2019,KhalilTudor2019b}, we consider linear equations driven by a space-time Gaussian noise with covariance structure generated by the Riesz kernel of order $4\gamma$ with $\gamma\in(0,1/4)$. Assuming the same sampling scheme as in the bounded domain case, we derive consistent and asymptotically normal estimators for $\theta$ or $\sigma$. We remark that the obtained results hold true for any $\alpha>0$, generalizing the results of \cite{KhalilTudor2019,KhalilTudor2019b}, where it is assumed that $\alpha\in(0,2]$. The case of nonlinear equations on the whole space is omitted in this study due to the lack of results on fine regularity properties of the solution (the so-called $L^p$ theory). We validate the theoretical results by numerical simulations for various sets of parameters; see Section~\ref{sec:examples}. In particular, we compute explicitly the aforementioned bias, which indeed turns out to be a significant correction to the naively derived  estimators.

\section{Preliminaries}\label{sec:prelim}

We fix a complete probability space $\bF=(\Omega, \sF, \bP)$ and throughout, all equalities and inequalities are understood in $\bP$-a.s. sense, unless otherwise stated.
As usual, we will denote by $\bP-\lim$ or $\xrightarrow{\mathbb{P}}$ the convergence in probability, and $\wlim$ or $\xrightarrow{d}$ will stand for the convergence in distribution. Correspondingly, $a_n=o_{\bP}(b_n)$ means that $a_n/b_n \xrightarrow{\mathbb{P}} 0$. Moreover, we write $a_n\lesssim b_n$, if there exists a constant $C$, independent of $n$, such that $a_n\leq Cb_n$ for all $n\in\bN$.

Let $X_t, t\in\bR$, be a real valued measurable function, and denote by $J$, and $\Delta_h$, the integral, and respectively the difference operators of the form 
\begin{align*}
JX_t&:=\int_0^t X_r \dif r,\quad t\in \mathbb{R},\\
\Delta_h X_t&:=X_{t+h}-X_t,\quad t\in \mathbb{R},\ h>0. 
\end{align*}
As usual, we put $J^0X:=X$, and for $m\in\bN$, we define $J^mX:=JJ^{m-1}X$. Similar notations apply to $\Delta_h$.  Note that, for $M,m\in\bN_0:=\bN\cup\set{0}$, we have 
\begin{align*}
\Delta^M_h (J^m X_t)=\sum_{k=0}^{M}(-1)^{M-k}{M\choose{k}} J^m X_{t+kh},\quad t\in \mathbb{R},\ h>0. 
\end{align*}

We will denote by $C(\bR)$ the space of continuous and bounded functions on $\bR$ endowed with $\sup$-norm $\norm{f}_\infty:=\sup|f|$. Correspondingly, for $k\in\bN$, we put
$C^k(\bR):=\set{f\in C(\bR) \ : \ \norm{f}_{C^k(\bR)}:= \sum_{j\leq k} \norm{D^jf}_\infty <\infty }$, where $D$ stands for differential operator. 

One of the key ideas of this paper is tracking and using the classical regularity of a continuous function in terms of conveniently chosen integral and difference operators.  For this purpose, we will be using the \textit{H\"older-Zygmund spaces} $\cC^s(\R), \, s>0$, introduced in \cite{Zygmund1945}  and endowed with the norm
$$
\lVert f\rVert_{s}^{(k, M)} = \lVert f\rVert_{C^k(\R)} + \lvert f\rvert_s^{(k,M)},
$$ 
with
\begin{equation}\label{eq:HS-norm1}
\lvert f\rvert_s^{(k, M)} = \sup_{h>0}h^{-(s-k)} \norm{\Delta_h^{M} D^kf}_\infty, 
\end{equation}
and where $k\in\N_0$, $M\in\bN$, such that  $k<s$ and $M>s-k$. It can be shown (cf. \cite[Section1.2.2]{Triebel1992}), that for any such $k$ and $M$, and fixed $s$ the norms $\lvert \,\cdot\,\rvert_s^{(k, M)}$ are equivalent.  We also recall that for any $s>0$, $\cC^s(\R)$ coincides with the Besov space $B^s_{\infty, \infty}(\R)$ (see also \cite{GineNickl2015}), and for $s\notin\N$,  $\cC^s(\R)$ coincide with the classical H\"older spaces.
Thus, the H\"older-Zygmund norms measure the regularity of a continuous function in the classical sense.  In this study, we will be mainly interested in the case $k=0$, which corresponds to statistical experiment of discrete measurements of the underlying process itself. However, if the observer evaluates discretely some derivative of $f$, then one should consider  $k\geq1$. Thus, we emphasize that the choice of $k=0$ is primarily driven by practical reasons, but in principle all results can be elevated to the general case $k\in\bN_0$.

\section{Smooth perturbations of higher order power variations}\label{sec:PowerVar}

Let $\pi=\set{t_0,\ldots,t_{N}}$ be the uniform partition of size $N$ of the interval $[a,b]\subset[0,T]$, and put $h:=h_N:=(b-a)/N = t_{k+1}-t_k, k =0,\ldots,N$. For fixed $s>0$, $q,M, N\in\bN$, such that $N>M$,  we define 
$$
V_{q,M,s,N}(X) := \frac{1}{b-a}\sum_{k=0}^{N-M} h
\left|\frac{\Delta_h^MX_{t_k}}{h^s}\right|^q.
$$
Similar to the power variation of a process, we are interested in the limiting behavior of $V_{q,M,s,N}$ as $N\to\infty$. 
The $\Delta$-power variation of order $(q,M,s)$ of process $X$ is defined as 
\begin{equation}\label{eq:p-varUni}
V_{q,M,s}(X):= \bP-\lim_{N\to\infty}V_{q,M,s,N}(X), 
\end{equation}
provided that the limit (in probability) exists.  Note that $V_{p,1,1}$ corresponds to the (normalized) power variation of order $p$.

We start with a simple, but important, result that links the path continuity of the process $X$ with its  generalized power variation. 

\begin{lemma}\label{lem:ZygmundToVariation}
	Let $q, M\in\bN$, $s>0$, such that $M>s$. If $X\in \cC^s([a,b])$, then $V_{q, M, s, N}(X)$ is uniformly bounded in $N$. 
\end{lemma}
\begin{proof}
	This follows at once by noticing that 
$$
V_{q,M,s,N}(X) = \frac{1}{b-a}\sum_{k=0}^{N-M}(t_{k+1}-t_k)\left|\frac{\Delta_h^MX_{t_k}}{(t_{k+1}-t_k)^s}\right|^q
\lesssim (h^{-s}\lVert\Delta_h^M X\rVert_\infty)^q
\lesssim (\lVert X\rVert_s^{(0,M)})^q.
$$
\end{proof}

We give the main results of this section, which in the nutshell says that the central limit theorems for $\Delta$-power variations of a stochastic process remain invariant under smooth perturbations; see also \cite[Proposition~2.1]{CialencoHuang2017}. 

\begin{theorem}\label{thm:PerturbationVariation}
Let $q\geq 1$, $s>0$, $M\in\N$ with $M>s$. Assume that $X\in \cC^s([a,b])$ and for some $\alpha>0$, $\Sigma\geq 0$, the following limit exists 
\begin{equation}\label{eq:VariationGeneralCLT_original}
	\lim_{N\to\infty } h_N^{-\alpha}\left(V_{q, M, s, N}(X)-V_{q, M, s}(X)\right)\overset{d}=\mathcal{N}(0, \Sigma), 
\end{equation}
where $\cN(0,\Sigma)$ is a Gaussian random variable with mean zero and variance\footnote{As usual, zero variance case is interpreted as the Dirac point mass at the mean.} $\Sigma$. 
Then,	for any  $Y\in C^{s+\eta}([a,b])$ with  $\eta>\alpha$, and $M>s+\alpha$, 
\begin{equation}\label{eq:VariationGeneralCLT_perturbed}
	\lim_{N\to\infty} h_N^{-\alpha}\left(V_{q, M, s, N}(X+Y)-V_{q, M, s}(X)\right)\overset{d}=\mathcal{N}(0, \Sigma). 
\end{equation}
\end{theorem}
\begin{proof} 
Without loss of generality, we assume that $M>s+\eta$, otherwise take $\eta'$ instead of $\eta$ with $\alpha<\eta'<\eta\wedge(M-s)$. We proceed analogous to \cite[Proposition~2.1]{CialencoHuang2017}. It suffices to show
	\begin{align}\label{eq:VariationGeneralCLTPerturbationToZero}
	\lim_{N\to \infty} h_N^{-\alpha}(V_{q, M, s, N}(X+Y)-V_{q, M, s, N}(X))=0, \quad \textrm{a.s.}.
	\end{align}
	Let $g_N(r)=\left(V_{q,M,s,N}(X)^{1/q}+rV_{q,M,s,N}(Y)^{1/q}\right)^q$. Then, by Minkowski's inequality, $g_N(-1)\leq V_{q,M,s,N}(X+Y)\leq g_N(1)$, and there exist $\xi_1, \xi_2\in [0,1]$ (dependent on $N\in\N$) such that
	$$
	g_N'(-\xi_1)=g_N(-1)-g_N(0)\leq V_{q,M,s,N}(X+Y)-V_{q,M,s,N}(X)\leq g_N(1)-g_N(0) = g_N'(\xi_2).
	$$
Thus, it remains to show $h_N^{-\alpha} \sup_{-1\leq r\leq 1}g_N'(r)\xrightarrow{a.s.} 0$, as $N\rightarrow\infty$. For $r\in [-1,1]$,
$$
\left|g_N'(r)\right|\leq q\left|V_{q,M,s,N}(X)^{1/q}+rV_{q,M,s,N}(Y)^{1/q}\right|^{q-1}V_{q,M,s,N}(Y)^{1/q}\lesssim h_N^{\eta}V_{q,M,s+\eta,N}(Y)^{1/q},
$$
and by Lemma~\ref{lem:ZygmundToVariation} and $Y\in \cC^{s+\eta}([a,b])$, $V_{q,M,s+\eta,N}(Y)$ is bounded uniformly in $N$. The claim follows from $\alpha < \eta$.  
\end{proof}

\begin{remark} (i) We note that the restriction $M>s+\alpha$ can be always satisfied by choosing $M$ large enough. (ii) If $\Sigma=0$, then the limits  \eqref{eq:VariationGeneralCLT_original} and \eqref{eq:VariationGeneralCLT_perturbed} can be equivalently understood as limits in probability. This in turn can be re-formulated in the terms of rates of convergence, as we do, for example, in  Theorems~\ref{thm:BoundedMixedVariationAsymptotics} and \ref{th:SPDE-parmEst}.  (iii) The results in this section can be easily extended to $\Delta$-power variations over arbitrary sequence of partitions, not necessarily uniform. Namely, one can replace the sequence of uniform partitions with a sequence of partitions with vanishing mesh-size in the above limits. However, generally speaking the counterpart of limit \eqref{eq:p-varUni} (if exists), may depend on the choice of the sequence of partitions. 

\end{remark}

\section{The case of fBM}\label{sec:fBM}

We start by recalling that a fractional Brownian motion (fBm) with Hurst index $H\in(0,1)$ is a centered Gaussian process $B^H=(B^H_t)_{t\in \mathbb{R}}$ such that
$$
\mathbb{E}\left(B_t^HB_r^H\right)=\frac{1}{2}\left(|t|^{2H}+|r|^{2H}-|t-r|^{2H}\right), 
\quad t,r\in \mathbb{R}.
$$
A continuous stochastic process $X$ is called \textit{$s$-self-similar} or self-similar of index $s$ (or self-similar for short) if the law of $(h^{-s}X_{ht})_{t\in\bR}$ on $C(\R)$ does not depend on $h>0$.  The process $X$ is said to be {\it stationary} if the law of $(X_{t+u})_{t\in\bR}$ on $C(\R)$ does not depend on $u\in\bR$, and $X$ is said to have {\it stationary increments} if $\Delta_hX$ is stationary for all $h>0$. A fractional Brownian motion $B^H$ is a prominent example of a self-similar process (of index $H$) with stationary increments. Many core properties of fBm are directly linked to these two features. 
However, generally speaking, differences of integrals of fBm are not self-similar in the usual sense, but rather, one has to account for the step-width of the difference operator. Towards this end, we extend the notion of self-similarity to parametrized family of processes, say $X^{(h)}$, $h>0$. Primarily, we will be interested in parametrized family of process of the form $X^{(h)}=\Delta_h^MY$, where $M\in\bN_0$ and $Y$ is a process that does not depend explicitly on $h>0$. We say that a parametrized family of process $X^{(h)}$ is {\it parametrized $s$-self-similar} (or just parametrized self-similar) if the law of $(h^{-s}X^{(h)}_{ht})_{t\in\bR}$ is independent of $h>0$. 
We also note that in general, if $X$ is stationary, then $JX$ is not necessarily stationary.

\begin{lemma}\label{lemma:PropertiesDeltaJ}
Let $X$ and $X^{(h)}, \ h>0$, be centered Gaussian processes. Then: 
		\begin{enumerate}[(i)]
			\item $\Delta_h^2JX = \Delta_hJ\Delta_hX$.
			\item If $X$ is $s$-self-similar, then $JX$ is $(s+1)$-self-similar.
			\item If $X^{(h)}$ is parametrized $s$-self-similar, then $\Delta_hX^{(h)}$ is parametrized $s$-self-similar and $JX^{(h)}$ is parametrized $(s+1)$-self-similar.
			\item If $X$ is stationary, then $\Delta_hX$ and $\Delta_hJX$ are stationary for any $h>0$.
		\end{enumerate}
\end{lemma}
	\begin{proof}
First we note that if $X$ is a centered Gaussian process, then $JX$ and $\Delta_hX$ are also Gaussian and centered. 	Thus,  the law of $\Delta_hJX^{(h)}$ is determined by  $\mathbb{E}\left[\Delta_hJX^{(h)}_t\Delta_hJX^{(h)}_r\right]$, which is equal to $\int_t^{t+h}\int_r^{r+h}\mathbb{E}\left[X^{(h)}_{v}X^{(h)}_w\right]\dif v\dif w$, $t,r\in\bR$. 
Using this, the above properties follow now by direct calculations. 
\end{proof}

Next, we state some properties specific to integro-differences of fBm of the $J^mB^H$ and $\Delta_h^MJ^mB^H$.

\begin{lemma}\label{lemma:PropertiesDeltaJBmH} 
The following assertions hold true: 
	\begin{enumerate}[(i)]
		\item For $m\in\bN_0$ and $t,r\in \mathbb{R}$, we have
		\begin{align}\label{eq:CovBmH}
		\mathbb{E}\left(J^mB^H_t\cdot J^mB^H_r\right)&=\sum_{k=0}^m \frac{(-1)^k\left(t^{m-k}r^{m+k+2H}+r^{m-k}t^{m+k+2H}\right)}{2(m-k)!\displaystyle\prod_{i=1}^{m+k} (2H+i)}+\frac{(-1)^{m+1}|t-r|^{2m+2H}}{2\displaystyle\prod_{i=1}^{2m} (2H+i)}.
		\end{align}
		In addition, $J^mB^H$ is $(m+H)$-self-similar. By convention, $\displaystyle\prod_{i=1}^0(2H+i)=1$.
		
		\item For $M,m\in\bN_0$ and $t, r\in\bR$, we have
		\begin{align}\label{eq:CovDeltaJBmH}
		\mathbb{E}\left[\Delta_h^MJ^mB^H_t\Delta_h^MJ^mB^H_r\right] 
		&=\sum_{k,l=0}^{M}(-1)^{2M-k-l}{M\choose{k}}{M\choose{l}} \mathbb{E}\left[J^mB^H_{t+kh}J^mB^H_{r+lh}\right].
		\end{align}
	\item If $M\geq m$, then $\Delta_h^MJ^mB^H$ is parametrized $(m+H)$-self-similar and has stationary increments.
	\end{enumerate}
\end{lemma}
\begin{proof} 
(i)  We prove \eqref{eq:CovBmH} by induction in $m$. For $m=0$, \eqref{eq:CovBmH} is immediate.
For $m=1$, by direct computations, we have
		\begin{align*}
		\mathbb{E}\left(JB_t^H\cdot JB_r^H\right)&=\int_0^t\int_0^r \mathbb{E}\left(B_u^HB_v^H\right)\dif u\dif v\\
		&=\frac{1}{2}\int_0^t\int_0^r \left(u^{2H}+v^{2H}-|u-v|^{2H}\right)\dif u \dif v\\
		&=\frac{1}{2}\left[\frac{t\cdot r^{2H+1}+r\cdot t^{2H+1}}{1!(2H+1)}+\frac{|t-r|^{2H+2}-r^{2H+2}-t^{2H+2}}{(2H+1)(2H+2)}\right],
		\end{align*}
and hence \eqref{eq:CovBmH} is true for $m=1$. 
		Suppose \eqref{eq:CovBmH} holds true for $m\geq 0$. Then,
		\begin{align*}
		&\mathbb{E}\left(J^{m+1}B^H_t\cdot J^{m+1}B^H_r\right)=\int_0^t\int_0^r \mathbb{E}\left(J^mB^H_u\cdot J^mB^H_v\right)\dif u\dif v\\
		&=\int_0^t\int_0^r \left[\sum_{k=0}^m \frac{(-1)^k\left(v^{m-k}u^{m+k+2H}+u^{m-k}v^{m+k+2H}\right)}{2(m-k)!\prod_{i=1}^{m+k} (2H+i)}+\frac{(-1)^{m+1}|v-u|^{2m+2H}}{2\prod_{i=1}^{2m} (2H+i)}\right]\dif u\dif v\\
		&=\sum_{k=0}^{m+1} \frac{(-1)^k\left(t^{m+1-k}r^{m+1+k+2H}+t^{m+1-k}r^{m+1+k+2H}\right)}{2(m+1-k)!\prod_{i=1}^{m+1+k} (2H+i)}+\frac{(-1)^{m+2}|t-r|^{2(m+1)+2H}}{2\prod_{i=1}^{2(m+1)} (2H+i)},
		\end{align*}
		and thus \eqref{eq:CovBmH} is proved. Consequently, $(m+H)$-self-similarity of $J^mB^H$ follows from Lemma~\ref{lemma:PropertiesDeltaJ}(ii).
		
\smallskip
\noindent		(ii) Identity \eqref{eq:CovDeltaJBmH} is immediate.

\smallskip
\noindent		(iii) The parametrized self-similarity follows from Lemma~\ref{lemma:PropertiesDeltaJ}(iii). 
Finally, Lemma~\ref{lemma:PropertiesDeltaJ}(iv) yields stationarity for $\Delta_h^{M+1}J^mB^H=\Delta_h^{M-m}(\Delta_hJ)^m\Delta_hB^H$, where we use Lemma~\ref{lemma:PropertiesDeltaJ}(i) and the fact that $\Delta_hB^H$ is stationary for all $h>0$. 
The proof is complete.
\end{proof}

Let us fix $M\in\bN$ and $s>0$, and write $s=m+H$ with $m\in\bN_0$ and $H\in (0,1)$. In view of Lemma~\ref{lemma:PropertiesDeltaJBmH}, there exists $\mu_{M,s}>0$ such that
\begin{align*}
\mathbb{E}\Big|\Delta_h^MJ^mB_t^H\Big|^2=\mu_{M,s} h^{2s},
\end{align*}
for all $t\in\bR$ and $h>0$, and where $\mu_{M,s}$ is given by
\begin{align*}
\mu_{M,s} &:=\sum_{k=0}^{M} {M\choose{k}}^2\sum_{p=0}^m \frac{(-1)^{p}k^{2s}}{(m-p)!\prod_{i=1}^{m+p}(2H+i)}\\
& \qquad +\sum_{0\leq j<k\leq M} (-1)^{2M-k-j}  {M\choose{k}}{M\choose{j}}\Bigg[\frac{(-1)^{m+1}(k-j)^{2s}}{\prod_{i=1}^{2m}(2H+i)}\\
& \qquad \qquad +\sum_{p=0}^m\frac{(-1)^{p}\left(k^{m-p}j^{m+p+2H}+j^{m-p}k^{m+p+2H}\right)}{(m-p)!\prod_{i=1}^{m+p}(2H+i)}\Bigg].
\end{align*}
We further set 
\begin{align}
\rho_{M,s}(\ell):=\mu^{-1}_{M,s} h^{-2s}\mathbb{E}\left(\Delta^{M}_hJ^mB^H_{t}\cdot \Delta^{M}_hJ^mB^H_{t+h\ell}\right),\quad \ell\in\bN_0. 
\end{align}
Note that due to self-similarity and  stationary increments property of fBM, we have that $\rho_{M,s}(\ell)$ does not depend on $t\in\bR$ and $h>0$.

Next, we will investigate the asymptotic behavior of the $q$-th (Hermite) variation of $\Delta_h^MJ^mB^H$, for which we will make use of (Breuer-Major) Theorem~\ref{th:Breuer-Major} applied to process 
$Y_t=\left(\mu_{M,s}^{1/2}h^s\right)^{-1}\Delta^M_hJ^mB^H_t$. First we note that by Lemma~\ref{lemma:PropertiesDeltaJBmH} the process $Y$ is a centered stationary Gaussian process with unit variance. Next result will be used to show that \eqref{eq:sumEll} is satisfied.

\begin{lemma}\label{lemma:RhoAsymptotics}
Assume that $M,q\in\bN$ and $0<s<M-\displaystyle\frac{1}{2q}$.  Then 
	\begin{align}\label{eq:AutocorrelationSummable}
		\sum_{\ell\in\Z}|\rho_{M,s}(\ell)|^q<\infty.
	\end{align}
\end{lemma}

\begin{proof}
Without loss of generality, we assume that $\ell\geq M$. The covariance function $\rho_{M,s}(\ell)$ becomes
\begin{align}
\rho_{M,s}(\ell)
&=\mu_{M,s}^{-1}\sum_{0\leq j,k\leq M} (-1)^{2M-k-j}{M\choose{k}}{M\choose{j}}\Bigg[\frac{(-1)^{m+1}(j+\ell-k)^{2m+2H}}{\prod_{i=1}^{2m}(2H+i)}\nonumber\\
&\qquad \qquad +\sum_{p=0}^m\frac{(-1)^p\left(k^{m-p}(j+\ell)^{m+p+2H}+(j+\ell)^{m-p}k^{m+p+2H}\right)}{(m-p)!\prod_{i=1}^{m+p}(2H+i)}\Bigg]\nonumber\\
&=:c_1\Delta_1^{2M}f_1(\ell) + \sum_{p=0}^m\left[c_{2,p}\Delta_1^Mf_{2,p}(\ell)+c_{3,p}\Delta_1^Mf_{3, p}(\ell)\right],
\end{align}
where 
$$
f_1(x)=(x-M)^{2m+2H}, \quad f_{2,p}(x)=x^{m+p+2H}, \quad f_{3,p}(x)=x^{m-p}. 
$$ 
First note that 
$$
c_{2, p}=((m-p)!\prod_{i=1}^{m+p}(2H+i))^{-1}\sum_{k=0}^M(-1)^{M-k}{M\choose k}k^{m-p}=0,
$$ 
where we used the fact that $M>m$. By direct computations, one can show that $c_1\neq 0$, and $c_{3,p}\neq 0$, for any $H\in(0,1)$. It is clear that, as $\ell\to\infty$,  
$$
\Delta_1^{2M}f_1(\ell)=\Delta_1^{2M}J^{2M}f_1^{(2M)}(\ell)=(\Delta_1J)^{2M}f_1^{(2M)}(\ell)\sim f_1^{(2M)}(\ell)\sim \ell^{2m+2H-2M}.
$$
If $M\leq m+p$, we similarly deduce that $\Delta_1^Mf_{3,p}(\ell)\sim \ell^{m-p-M}$, and $\Delta_1^{2M}f_1$ grows faster than $\Delta_1^Mf_{3, p}$, since $2m+2H-2M>m-p-M$. If $M>m+p$, we have $\Delta_1^Mf_{3,p}(\ell)\equiv 0$. Combining the above, we have 
$$
\rho_{M,s}(\ell)\sim\ell^{2m+2H-2M}, \quad \ell\rightarrow\infty.
$$ 
Thus, if $H<M-m-\displaystyle\frac{1}{2q}$, then \eqref{eq:AutocorrelationSummable} is true. This concludes the proof.  
\end{proof}

As an immediate consequence of Lemma~\ref{lemma:RhoAsymptotics}, we get that for any $0<s<M-\displaystyle\frac{1}{2q}$, the quantity
\begin{align}
\rho^2_{q,M,s}:=q!\displaystyle\sum_{\ell\in \mathbb{Z}}\left(\rho_{M,s}(\ell)\right)^q
\end{align}
is well-defined and finite.

The following result identifies $V_{q,M,s}(J^mB^H)$ for $s=m+H$ together with its convergence rate.  

\begin{theorem}\label{thm:itfBmCLT}
Let $M>m\geq 0$ and $q\geq 1$ be integers, and assume that either of the following assumptions is satisfied:
	\begin{enumerate}[(i)]
		\item $M=m+1$ and $0<H<3/4$, 
		\item $M\geq m+2$ and $0<H<1$. 
	\end{enumerate}
Then,  there exists $\sigma_{q,M,s}>0$ such that 
	\begin{align}\label{eq:itfBmCLT-1}
		\sqrt{N}\left(V_{q,M,s,N}\left(J^mB^H\right)-\tau_q\mu_{M,s}^{q/2}\right)\xrightarrow{d}\mathcal{N}\left(0,\sigma_{q,M,s}^2\mu_{M,s}^q\right), \quad \textrm{as} \ N\to\infty,
	\end{align}
where  $\tau_q:=\mathbb{E}|Z|^q$ with $Z\sim\mathcal{N}(0,1)$.

Moreover, if $q$ is an even number, then $\sigma_{q,M,s}^2=\sum_{k=1}^q {q\choose{k}}^2\tau_{q-k}^2\rho^2_{k,M,s}$.
\end{theorem}
\begin{proof} We apply  Theorem~\ref{thm:Breuer-Major}, by taking  $(Y_k)_{k\in\Z}=\left(\mu_{M,s}^{-1/2}h^{-s}\Delta_h^MJ^mB^H_{t_k}\right)_{k\in\Z}$ and  $f(x)=|x|^q-\tau_q=\sum_{k=0}^\infty a_kH_k(x)$ with $a_k=(2\pi)^{-1/2}\int (|x|^q-\tau_q)H_k(x)e^{-x^2/2}\dif x$. Note that, in view of \cite[Example 7.2.2]{NourdinPeccati2012} the function $f$ has Hermite rank $d=2$, namely $a_0=a_1=0$ and $a_2\neq 0$. It remains to show that \eqref{eq:sumEll} is satisfied, which in our case becomes $\sum_{\ell\in\Z}\rho^2_{M,s}(\ell)<\infty$. 
By Lemma~\ref{lemma:RhoAsymptotics}, this is true if $0<s<M-1/4$, or equivalently if $0<H<(M-m)-1/4$, which is  satisfied in view of assumptions (i)-(ii).  Thus, \eqref{eq:itfBmCLT-1} is proved. 

For $q$ even, it can be shown (for example, by induction, or see \cite[p.1076]{NourdinNualartTudor2010}) that 
	\begin{align}
	V_{q,M,s,N}\left(J^mB^H\right)-\frac{N-M+1}{N}\tau_q\mu_{M,s}^{q/2}=\frac{1}{N}\mu_{M,s}^{q/2}\sum_{k=1}^q{q\choose {k}} \tau_{q-k}\sum_{j=0}^{N-M} H_k\left(\frac{\Delta_h^MJ^mB_{t_{j}}^H}{\mu_{M,s}^{1/2}h^{s}}\right).
	\end{align}
\end{proof}

\begin{remark} (i) 	We emphasize that the limit of $V_{q, M, s, N}(J^mB^H)$ depends through $\mu_{M,s}$ on the regularity $s$ of the process as well as the number of differences $M$. In particular, even for small $h$ it is not possible to approximate the rescaled finite difference operator $h^{-1}\Delta_h$ by a derivative operator without introducing a non--trivial bias. 
(ii) The constant $\mu_{M,s}$ can be easily computed, and for reader's convenience we list some of its values. If $M=1$, $m=0$ and $0<H<3/4$, then $\mu_{M,s}=1$. If $M=2$, $m=1$ and $H=1/4$, then $\mu_{M,s}=(\sqrt{2}-1)\frac{16}{15}\approx 0.44$. If $M=2$, $m=1$ and $H=1/2$, then $\mu_{M,s}=2/3$.

\end{remark}

\section{Semilinear SPDEs on a bounded domain}\label{sec:SPDE-BDom}
In this section we consider SPDEs on $\mathcal{D}=(0,1)$ with zero boundary conditions. Towards this end, for $k\geq 1$, set $\Phi_k(x)=\sqrt{2}\sin(k\pi x)$ and $\lambda_k=k^2\pi^2$. The set  $\set{\Phi_k}_{k\in\bN}$ forms an orthonormal basis in $L^2(\mathcal{D})$. Further, for $s\in\R$, set $H^s(\mathcal{D}):=\{u\in L^2\;|\;\sum_{k=1}^\infty \lambda_k^{s}(u,\Phi_k)^2<\infty\}$. 
The Laplacian $\Delta = \partial_{xx}$, acting on $C^\infty(\mathcal{D})$, can be extended to a closed, densely defined operator $\Delta$ on $L^2(\mathcal{D})$ with domain $H^2(\mathcal{D})$ and compact resolvent. The $\Phi_k$ are eigenfunctions of $-\Delta$ with eigenvalues $\lambda_k$. Not that for $s>1/2$, $H^s(\mathcal{D})\rightarrow C(\mathcal{D})$, and $u(0)=u(1)=0$ for $u\in H^s(\mathcal{D})$. 

We consider the following semilinear SPDE on $L^2(\mathcal{D})$:
\begin{equation}\label{eq:SPDE-BDom}
	\mathrm{d}X_t = \left(-\theta(-\Delta)^{\alpha/2} X_t + F(X_t)\right)\mathrm{d}t + \sigma B\mathrm{d}W_t, \quad X_0\in L^2(\cD),
\end{equation}
where $\alpha,\theta, \sigma > 0$, $W$ is a cylindrical Wiener process on $L^2(\mathcal{D})$, $B=(-\Delta)^{-\gamma}$ for some $\gamma>1/4$, and $F$ is a nonlinear operator.

We assume that \eqref{eq:SPDE-BDom} is well-posed, and refer, for instance, to \cite{DaPratoZabczykBook2014,LiuRoeckner2015} for sufficient conditions.  
The condition $\gamma>1/4$ is imposed to have function valued solutions, which in turn is used to define and interpret naturally the nonlinear part $F(X)$; see Example~\ref{exampleSPDE}. In principle, for some classes of equation, such as linear equations, one can consider less restrictive values for $\gamma$ and take, for instance, $\gamma=0$ that will correspond to space-time white noise. Since our focus is mainly on nonlinear equations, we omit discussing these cases herein. 

As customary in statistical inference for nonlinear SPDEs \cite{IgorNathanAditiveNS2010,PasemannStannat2019,ACP2020}, we will use the splitting of the solution argument, by writing $X=\overline X + \widetilde X$, where
\begin{align}
	\mathrm{d}\overline X_t &= -\theta (-\Delta)^{\alpha/2}\overline X_t\mathrm{d}t + \sigma B\mathrm{d}W_t,\quad \overline X_0=0, \label{eq:DirichletConditionLinearEquation}\\
	\mathrm{d}\widetilde X_t &=\left(-\theta (-\Delta)^{\alpha/2}\widetilde X_t + F(\overline X_t + \widetilde X_t)\right)\mathrm{d}t,\quad \widetilde X_0=X_0. \label{eq:barXBoundDom}
\end{align}

The solution to \eqref{eq:DirichletConditionLinearEquation} can be expressed either as a Fourier series, or can be given by the stochastic convolution
\begin{equation}\label{eq:SPDE-Fourier}
\overline X_t = \sigma\int_0^t e^{-\theta(t-r)(-\Delta)^{\alpha/2}}B\mathrm{d}W_r = \sum_{k=1}^\infty\left(\sigma\lambda_k^{-\gamma}\int_0^te^{-\theta(t-r)\lambda_k^{\alpha/2}}\mathrm{d}W^{(k)}_r\right)\Phi_k=:\sum_{k=1}^\infty \overline x_k(t)\Phi_k,
\end{equation}
where $W^{(k)}=(W,\Phi_k)_{L^2}, k\geq 1$, are independent one-dimensional Brownian motions, $t\mapsto e^{-\theta t(-\Delta)^{\alpha/2}}$, $t>0$, is the $C_0$-semigroup on $L^2(\mathcal{D})$ generated by $-\theta(-\Delta)^{\alpha/2}$, and the convergence is understood in $L^2(\mathcal{D})$. Note that $\overline{x}_k(t)=(\overline{X}(t), \Phi_k)_{L^2}$, i.e. $\overline{x}_k(t)$ is also the Fourier coefficient of the solution $\overline{X}(t)$ with respect to $\set{\Phi_k}_{k\in\bN}$. 

The next two results provide some fine continuity properties of the trajectories of $\overline{X}$ and $\widetilde{X}$.

\begin{proposition}\label{prop:RegXbar}
For any $s<2\gamma+\alpha/2-1/2$, it holds that $\overline X\in C(0, T; C^s(\mathcal{D}))$.
\end{proposition}
\begin{proof} The common line of attack is to show that $\overline{X}\in C(0,T;W^{s, p}(\mathcal{D}))$ for any $p\geq2$, and then employ the Sobolev embedding theorem.  We refer, for example, to \cite[Appendix B.1]{ACP2020} for details when $\alpha=2$, and since the proof is based on the Fourier decomposition of the solution in the base $\set{\Phi_k}$, the general case is obtained similarly. 
\end{proof}

As a direct consequence, we note that  $\overline X$  has up to $\lfloor 2\gamma+\alpha/2-1/2\rfloor$ classical derivatives. We call $s^*=2\gamma+\alpha/2-1/2$ the optimal regularity, and we make a  standing assumption that $s^*>0$ and $s^*\notin \bN$. 

\begin{proposition}\label{prop:HigherNonlinearRegularity}
	Assume that there exist $\eta, \epsilon>0$, $0\leq s_0<s^*$, and a continuous function $g:[0,\infty)\rightarrow[0,\infty)$, such that for any $s_0\leq s < s^*$,  	
	\begin{equation}\label{eq:ConditionF}
		\lVert F(u)\rVert_{s+\eta {-\alpha}+\epsilon}\leq g(\lVert u\rVert_{s}),
	\end{equation}
	where, as before, $\lVert\cdot\rVert_s$ denotes the H\"older-Zygmund norm. 
	Then,  $\widetilde X\in C(0,T;C^{s+\eta}(\mathcal{D}))$, for any $0\leq s < s^*$.
\end{proposition}

\begin{proof} The case of Sobolev spaces $W^{s,p}(\mathcal{D})$ instead of H\"older spaces $C^s(\mathcal{D})$ has been treated in  \cite{ACP2020,PasemannEtAl2020} for $\alpha=2$, the proof for general $\alpha$ is identical. Since for an arbitrary chosen $\epsilon>0$, we have that $C^s(\mathcal{D})\subset W^{s, p}(\mathcal{D})\subset C^{s-\epsilon}(\mathcal{D})$ for large enough $p$, the desired result follows at once. 
\end{proof}

\begin{example}\label{exampleSPDE}
We present several types of semilinear SPDEs whose nonlinearity $F$ satisfies \eqref{eq:ConditionF}, which in particular guarantees that all results from this section hold for the solutions to these classes of equations.  For technical details see \cite{PasemannStannat2019,ACP2020}. 
	\begin{enumerate}[1)]
		\item {\it (fractional) Heat equation:} In the case $F=0$, \eqref{eq:SPDE-BDom} becomes linear, sometimes called fractional head equation, and \eqref{eq:ConditionF} is trivially satisfied for any $\eta>0$. 
		\item {\it Reaction-diffusion equation:} Let $F(u)(x)=f(u(x))$, where $f$ is a polynomial function or $f\in C_b^\infty(\mathcal{\R})$. Then \eqref{eq:ConditionF} is true for any $0<\eta<2$. 
		\item {\it Advection-diffusion equation:} Let $F(u)=v\partial_xu$ for a given $v\in C^\infty(\mathcal{D})$. Then \eqref{eq:ConditionF} holds with any $0<\eta<1$. 
		\item If $F=F_1+F_2$, for some $F_1,F_2$ that satisfy \eqref{eq:ConditionF} with continuous functions $g_1, g_2$, then $F$ satisfies \eqref{eq:ConditionF} with $g=g_1+g_2$.
	\end{enumerate}
\end{example}

Next, using representation \eqref{eq:SPDE-Fourier}, we set 
$$
\xi_k:= \sqrt{\displaystyle\frac{2\theta \lambda_k^{\alpha/2+2\gamma}}{\left(1-e^{-2\theta \lambda_k^{\alpha/2} t}\right)\sigma^2}}\cdot \overline x_k(t), \qquad  k\geq 1. 
$$
Clearly, $\xi_k$'s are independent standard Gaussian random variables, and $\overline X_t$ can be written as
\begin{equation}\label{eq:SPDE-fBMdec}
\overline X_t=\frac{\sigma}{\sqrt{2\theta}}\sum_{k\geq 1}\frac{1}{\lambda_k^{\alpha/4+\gamma}}\xi_k\Phi_k+\frac{\sigma}{\sqrt{2\theta}}\sum_{k\geq 1}\frac{\left(\sqrt{1-e^{-2\theta \lambda_k^{\alpha/2}t}}-1\right)}{\lambda_k^{\alpha/4+\gamma}}\xi_k\Phi_k.
\end{equation}
Recall that under our standing assumption $s^*\notin\N_0$, and thus us write $s^*=m+H$ for some unique $m\in\N_0$ and $0<H<1$. As one may expect, $H$ will be linked to the Hurst parameter of a fBM. 
In particular, if $\alpha=2$ and $\gamma=0$, then $s^*=1/2$, $m=0$ and $H=1/2$, and as shown in~\cite{CialencoHuang2017} the first term in \eqref{eq:SPDE-fBMdec} is a fBM with Hurst index $H$ and the second term is an infinitely smooth process.

In view of \eqref{eq:SPDE-fBMdec}, we have, for fixed $t>0$,

\begin{equation}\label{eq:Deru-L1}
\partial_x^m\overline X_t  =
\begin{cases}
(-1)^{m/2} \frac{\sigma}{\sqrt{2\theta}} L_{\sin}^H+R_{\sin}, & \textrm{ if } m \textrm{ \ is even}, \\
(-1)^{(m-1)/2} \frac{\sigma}{\sqrt{2\theta}} L_{\cos}^H+R_{\cos}, & \textrm{ if } m \textrm{ \  is odd}, 
\end{cases}
\end{equation}
where 
\begin{align}
	L_{\sin}^H(x):=\displaystyle\sqrt{2}\sum_{k\geq 1}\lambda_k^{-H/2-1/4}\xi_k\sin(k\pi x), \quad L_{\cos}^H(x):=\displaystyle\sqrt{2}\sum_{k\geq 1}\lambda_k^{-H/2-1/4}\xi_k\cos(k\pi x), 
\end{align}
and $R_{\sin},R_{\cos}\in C^{\infty}\left(\mathcal{D}\right)$. Note that here we used that 
		\begin{equation}\label{eq:der-sum}
		\partial_x^m\sum_{k\geq 1}\frac{1}{\lambda_k^{\alpha/4+\gamma}}\xi_k\Phi_k(x)=\sum_{k\geq 1}\frac{1}{\lambda_k^{\alpha/4+\gamma}}\xi_k\partial_x^m\Phi_k(x),
		\end{equation} 
which is true thanks to the uniform converges of the last series. The latter is due to the following estimates 
\begin{align*}
\mathbb{E}\left|\sum_{k\geq 1}\frac{1}{\lambda_k^{\alpha/4+\gamma}}\xi_k\partial_x^m\Phi_k(x)\right|^2\lesssim \sum_{k\geq 1} k^{2m-\alpha-4\gamma}<\infty. 
\end{align*}

Motivated by \cite{Picard2011}, next we show that the stochastic processes $L_{\sin}^H$ and $L_{\cos}^H$, as functions of $x$, are strongly related to a fBM. 
For $0<H<1$, $H\neq 1/2$, let
\begin{align}
	\unitvar_H := -\frac{2}{\pi}\Gamma(-2H)\cos(\pi H),
\end{align}
and further put $\unitvar_H = 1$ for $H=1/2$. The constant $\nu_H$ corresponds to $\rho_H$ in \cite{Picard2011}.  We emphasis that the fBm in this work, say $B^H_x, x\geq0$, is scaled as in most of the literature, namely  $\bE[(B_1^H)^2]=1$, in contrast to $\cite{Picard2011}$, where the fBm is scaled such that $\bE[(B_1^H)^2]=\nu_H$. Respectively, some of the results from \cite{Picard2011} used below have to be adjusted accordingly.

\begin{lemma}\label{lem:DirichletBasisFBMRepresentation}
Let $0<H<1$, and $B^H_x, \ x\geq 0$, be a fBM with Hurst parameter $H$. There exists a stochastic process $R^H\in C^\infty(\R)$, such that for any $0<a<b<1$, the following hold true: 
	\begin{enumerate}[1)]
		\item The probability laws  of $\unitvar_H^{1/2}B^H_\cdot$ and $(L_{\sin}^H(a+\cdot) + R^H(a+\cdot))-(L_{\sin}^H(a)+R^H(a))$ are equivalent on (canonical space) $C([0, b-a])$.
		\item The laws of $\unitvar_H^{1/2}B^H_\cdot$ and $L_{\cos}^H(a+\cdot) - L_{\cos}^H(a)$ are equivalent on $C([0, b-a])$.
	\end{enumerate}
Moreover, if $H=1/2$, then above laws are even equal. 
\end{lemma}
\begin{proof}
Same as in \cite{Picard2011}, we define the process 
$$
\widehat{B}_x^H = \xi_0x + \sqrt{2}\sum_{k\geq 1} \left( \xi_k'  \frac{\cos(2\pi kx)-1}{(2\pi k)^{H+1/2}} + \xi_k'' \frac{\sin(2\pi kx)}{(2\pi k)^{H+1/2}}\right), 
$$
where $\xi_k',\xi_k''$ are i.i.d. standard normal random variables, and $x\in\bR$. Note that in view of \cite[Theorem~27]{Picard2011},  $\widehat B^H$ has stationary increments. Then, (i) is proved by following similar steps as in the proof of \cite[Theorem~30]{Picard2011}, with $\bar B^H$ replaced by $\widehat B^H$, and noting that 
\[
2^{H-1/2}\left(\widehat B^H_{x/2}-\widehat B^H_{-x/2}\right) = 2^{H-1/2}\xi_0x+L_{\sin}^H(x). 
\]
Consequently, taking $R^H(x)=2^{H-1/2}\xi_0 x$ we have proved (i). \\
(ii) We proceed similarly, and note that 
$$
2^{H-1/2}\left(\widehat B^H_{x/2}+\widehat B^H_{-x/2}\right) = L_{\cos}^H(x) - \frac{1}{\sqrt{2}}\sum_{k=1}^\infty\xi_k\lambda_k^{-H/2-1/4}=:L_{\cos}^H(x) - c^H.
$$
From here, since clearly $c^H\in L^2(\Omega)$, for any $H>0$, we conclude that the law of $L_{\cos}^H-c^H$ on $C([0,b])$ is equivalent to the law of $(\unitvar_H^{1/2}B^H_x+\unitvar_H^{1/2}B^H_{-x})/\sqrt{2}$. In view of \cite[Remark~5.11]{Picard2011}, the increments of this process, and thus also the increments of $L_{\cos}^H$, have a law equivalent to the law of $\unitvar_H^{1/2}B^H_{\cdot-a}$ on $C([a, b])$, and the claim follows. The case $H=1/2$ is known, and follows, for example, by Karhunen-Loeve type expansions of Brownian motion. 
\end{proof}

\begin{proposition}\label{prop:RepresentationX}
Let $m\in\N_0$ and $0<H<1$ such that $m+H=s^*=2\gamma+\alpha/2-1/2$. Then, there exists a stochastic process $R^{m,H}\in C^\infty(\mathcal{D})$, such that for any $0<a<b<1$, the laws of $(-1)^{\lfloor m/2\rfloor}\sigma^{-1}\unitvar_H^{-1/2}\sqrt{2\theta}\ \overline X+R^{m,H}$ and $J^mB^H_{\cdot-a}$ are equivalent on $C([a,b])$. If $H=1/2$, the laws are even equal. 
\end{proposition}
\begin{proof}
Applying $J^m$ to \eqref{eq:Deru-L1}, we note that it suffices to prove that for any $\bar m\in\N_0$ there exist $R_{\sin}^{\bar m,H}, R_{\cos}^{\bar m,H}\in C^\infty(\R)$ such that the laws $P_{\sin}^{(\bar m)}$ of $J^{\bar m}L_{\sin}^H+R_{\sin}^{\bar m, H}$ and $P_{\cos}^{(\bar m)}$ of $J^{\bar m}L_{\cos}^H+R_{\cos}^{\bar m, H}$ are equivalent to the law $Q^{(\bar m)}$ of $\unitvar_H^{1/2}J^{\bar m}B^H_{\cdot-a}$ on $C([a, b])$. 

We will prove the above by induction in $\bar m$. First, for $z\in\R$, let $\tau_z$ be the shift operator, given by $\tau_zf:=f_{\cdot+z}$, and we view $\tilde J:=\tau_{-a}J\tau_{a}$ as a bounded operator $\tilde J:C([a, b])\rightarrow C([a, b])$.   For $\bar m=0$, we note that by Lemma~\ref{lem:DirichletBasisFBMRepresentation}, $P_{\sin}^{(0)}$, $P_{\cos}^{(0)}$ and $Q^{(0)}$ on $C([a, b])$ are equivalent, where $R_{\sin}^{0,H}(x)=R^H(x)-(L_{\sin}^H(a)+R^H(a))$, $R_{\cos}^{0,H}(x)\equiv-L_{\cos}^H(a)$.  
Now assume that the claim is true for $\bar m\geq 0$. Then the pushforward measures $\tilde J^*P_{\sin}^{(\bar m)}$, $\tilde J^*P_{\cos}^{(\bar m)}$ and $\tilde J^*Q^{(\bar m)}$ are equivalent measures on $C([a, b])$. As $\tilde J=\tau_{-a}J\tau_a$, we see that $\tilde J^*Q^{(\bar m)}$ is the law of $\unitvar_H^{1/2}J^{\bar m}B^H_{\cdot-a}$, i.e. $\tilde J^*Q^{(\bar m)}=Q^{(\bar m+1)}$. Likewise, $\tilde J^*P_{\sin}^{(\bar m)}$ is the law of $J^{\bar m+1}L_{\sin}^H+R_{\sin}^{\bar m+1, H}$ with $R_{\sin}^{\bar m+1, H}=\tilde JR_{\sin}^{\bar m,H}-\int_0^aJ^{\bar m}L_{\sin}^H(y)\mathrm{d}y$. For this choice of $R_{\sin}^{\bar m, H}$ it holds $\tilde J^*P_{\sin}^{(\bar m)}=P_{\sin}^{(\bar m+1)}$. The case of $\tilde J^*P_{\cos}^{(\bar m)}$ is treated similarly. If $H=1/2$, one can trace the above arguments and notice that equivalent laws can be replaced with equal laws.  The proof is complete. 
\end{proof}

Now, we are in the position to prove the main result of this section. In the sequel, we fix $0<a<b<1$ and consider the generalized variation $V_{q,M,s,N}(X_t)$ on $[a,b]$, namely on an interval away from the boundary.  

\begin{theorem}\label{thm:BoundedMixedVariationAsymptotics}
	Let $M, q\in\N$, and assume that either $M=m+1$ with $H<1/2$ or $M\geq m+2$. 
	Suppose that \eqref{eq:ConditionF} holds for some $\eta> 1/2$. 
	Then, for any $\epsilon>0$,
	\begin{equation}\label{eq:DirichletTheoremVariationAsymptotics}
		V_{q,M,s^*,N}(X_t) = \tau_q\left(\frac{\sigma^2\unitvar_H\mu_{M, s^*}}{2\theta}\right)^{q/2} + o_\mathbb{P}(N^{-1/2+\epsilon}).
	\end{equation}
If in addition $s^*\in 1/2+\N_0$, then 
	\begin{align}\label{eq:DirichletTheoremVariationAsymptoticsCLT}
		\sqrt{N}\left(V_{q,M,s^*,N}(X_t) - \tau_q\left(\frac{\sigma^2\unitvar_H\mu_{M, s^*}}{2\theta}\right)^{q/2}\right)\xrightarrow{d}\mathcal{N}\left(0, \left(\frac{\sigma^2\unitvar_H\mu_{M, s^*}}{2\theta}\right)^q\sigma_{q,M,s^*}^2\right).
	\end{align}
	
\end{theorem}
\begin{proof}
Set $Z^{m,H}:=(-1)^{\lfloor m/2\rfloor}\sigma^{-1}{\unitvar_H^{-1/2}}\sqrt{2\theta}\overline X + R^{m, H}$, with $R^{m,H}$ as in Proposition~\ref{prop:RepresentationX}.
Since $\eta> 1/2$, Proposition~\ref{prop:HigherNonlinearRegularity} and Theorem~\ref{thm:PerturbationVariation} (with $\alpha=1/2-\epsilon$, $\Sigma=0$ or $\alpha=1/2$, {$\Sigma=\sigma^{2q}\unitvar_H^q\mu_{M,s^*}^q\sigma_{q,M,s^*}^2/(2\theta)^q$} in the notation therein) imply that it is enough to show that \eqref{eq:DirichletTheoremVariationAsymptotics}, and \eqref{eq:DirichletTheoremVariationAsymptoticsCLT} hold with $X_t$ replaced by $\overline X_t$. Consequently, since $R^{m,H}\in C^\infty(\mathcal{D})$, the claims are, respectively, equivalent to 
	\begin{align}
		N^{1/2-\epsilon}\left(V_{q,M,s^*,N}(Z^{m,H}) - \tau_q\mu_{M, s^*}^{q/2}\right)&\xrightarrow{\mathbb{P}}0, \label{eq:ThmDirichletRepresentationX} \\
		\sqrt{N}\left(V_{q,M,s^*,N}(Z^{m,H})-\tau_q\mu_{M,s^*}^{q/2}\right)&\xrightarrow{d}\mathcal{N}\left(0, \sigma_{q,M,s^*}^2\mu_{M,s^*}^q\right),\label{eq:ThmDirichletRepresentationXCLT}
	\end{align}
for a general $s^*$, and, respectively, for $s^*\in 1/2+\N_0$. From Theorem \ref{thm:itfBmCLT} it follows that 
\begin{equation}
\int_{C([a, b])}\mathbbm{1}\left(\left| N^{1/2-\epsilon}{\left(V_{q,M,s^*,N}(f) - \tau_q\mu_{M, s^*}^{q/2}\right)}\right|>\epsilon'\right)\mathrm{d}\mathcal{L}(J^mB^H_{\cdot-a})(f)\rightarrow 0,
\end{equation}
for any $\epsilon, \epsilon'>0$. 
By Proposition~\ref{prop:RepresentationX}, the laws of $J^mB^H_{\cdot-a}$ and $Z^{m,H}$ are equivalent on $C([a, b])$, and thus\footnote{We recall that for two equivalent measures $P\sim Q$ on some measureable space $(M, \mathcal{M})$, it holds that $P(A_N)\rightarrow 0$ if and only if $Q(A_N)\rightarrow 0$ for any $(A_N)_{N\in\N}\subset\mathcal{M}$, see e.g. \cite[Chapter 6]{Vaart1998}.}
	\begin{align}\label{eq:MixedConditionsConvergenceLaw}
		\int_{C([a, b])}\mathbbm{1}\left(\left| N^{1/2-\epsilon}{\left(V_{q,M,s^*,N}(f) - \tau_q\mu_{M, s^*}^{q/2}\right)}\right|>\epsilon'\right)\mathrm{d}\mathcal{L}(Z^{m, H})(f)\rightarrow 0
	\end{align}
	for any $\epsilon, \epsilon'>0$, which is equivalent to \eqref{eq:ThmDirichletRepresentationX}. Finally, if $s^*\in 1/2+\N_0$, then $Z^{m,H}$ and $J^mB^H_{\cdot - a}$ are equal in law, and \eqref{eq:ThmDirichletRepresentationXCLT} becomes \eqref{eq:itfBmCLT-1}. 
	This concludes the proof. 
\end{proof}

As a direct consequence, we obtain  a procedure to estimate one of the parameters $\sigma, \theta$, if the other one is known, based on discrete observations on the uniform grid of $[a, b]$. 

\begin{theorem}\label{th:SPDE-parmEst}
In the setting of Theorem~\ref{thm:BoundedMixedVariationAsymptotics}, the following hold true: 
	\begin{enumerate}[(i)]
		\item If $\theta$ is known, then $\widehat\sigma^{q, M}_N:=\tau_q^{-1}(2\theta/{(\unitvar_H\mu_{M, s^*})})^{q/2}V_{q,M,s^*,N}(X_t)$
		is a consistent estimator for $\sigma^q$, and for any $\epsilon>0$, 
		$$
		\widehat\sigma^{q,M}_N=\sigma^q + o_\mathbb{P}(N^{-1/2+\epsilon}).
		$$
If $s^*\in 1/2+\N_0$, then also
		\begin{align*}
			\sqrt{N}\left(\widehat\sigma^{q, M}_N - \sigma^q\right)\xrightarrow{d}\mathcal{N}\left(0, {\frac{\sigma^{2q}}{\tau_q^2}\sigma_{q,M,s^*}^2}\right).
		\end{align*}
		
		\item If $\sigma$ is known, then $\widehat\theta^{q, M}_N:=\tau_q^{2/q}{\unitvar_H}\mu_{M, s^*}\sigma^2 / (2V_{q, M, s^*, N}(X_t)^{2/q})$ is a consistent estimator for $\theta$, and 
		$$
		\widehat\theta^{q, M}_N = \theta + o_\mathbb{P}(N^{-1/2+\epsilon}),
		$$ for any $\epsilon>0$. If $s^*\in 1/2+\N_0$, then
		\begin{align*}
			\sqrt{N}\left(\widehat\theta^{q, M}_N - \theta\right)\xrightarrow{d}\mathcal{N}\left(0, {\frac{4\theta^2}{q^2\tau_q^2}\sigma_{q, M, s^*}^2}\right).
		\end{align*}
		
	\end{enumerate}
\end{theorem}
We conclude this section with several remarks: 
\begin{enumerate}
\item The choice of the Dirichlet boundary conditions is not essential. By changing the role of $L_{\sin}^H$ and $L_{\cos}^H$, we immediately get an analogous result for Neumann boundary conditions. Similarly, using the representation  of $B^H$  in terms of $L_{\mathrm{mix}}^H=\sum_{k\geq 1}{(\lambda_k^M)}^{-H/2-1/4}\xi_k\Phi_k^M$ (cf. \cite[Theorem 6.19]{Picard2011}) with $\Phi_k^M(x)=\sqrt{2}\sin((k-1/2)\pi x)$ and $\lambda_k^M=(k-1/2)^2\pi^2$, we get the same result for mixed boundary conditions. 

\item Applying the central limit theorem from Theorem~\ref{thm:itfBmCLT} to SPDEs, essentially depends on establishing a stronger than equivalence in law representation of the solution $X_t$ in terms of a fractional Brownian motion. To the best of our knowledge, this is an open problem for a general $H\neq 1/2$.

\item Similar results can be derived if \eqref{eq:SPDE-BDom} is driven by an additive space-only noise (the so-called parabolic Anderson model) instead of space-time noise. The main difference in this case is that the optimal regularity is $s^*=2\gamma+\alpha-1/2$ (cf. \cite{CialencoKimLototsky2019, CialencoKim2020}).  

\end{enumerate}

\section{Linear SPDEs on unbounded domain}\label{sec:unboundedDom}

We consider the (linear) counterpart of \eqref{eq:SPDE-BDom} on the whole space, namely the stochastic evolution equation of the form
\begin{align}\label{eq:fracHeatR}
	\partial_tX_t(x)&=-\theta \left(-\Delta\right)^{\alpha/2}X_t(x)+\sigma \dot W^{\gamma}(t,x),\quad t>0,\ x\in\mathbb{R},\\
	X_0(x)&=0,\quad x\in \mathbb{R},\nonumber
\end{align}
where $\alpha,\theta,\sigma>0$ and $W^\gamma(t,A), \ t\geq 0, A\in\cB(\bR)$, for some $\gamma\in(0,1/4)$, is a centered Gaussian field with covariance structure 
$$
\bE[W^\gamma(t,A) W^\gamma(s,B)] = (t\wedge s) \int_A\int_B K_\gamma(x-y)\dif x \dif y, 
$$
with $K_\gamma$ being the so-called Riesz kernel of order $\gamma$ given by 
$$
K_\gamma(x) =  \frac{\Gamma(1/2-2\gamma)}{2\pi^{3/2} \Gamma(2\gamma)} \cdot |x|^{4\gamma-1}, \quad \gamma\in(0,1/4). 
$$
We remark that traditionally in the literature the Riesz kernel has slightly different parameterization, with $\gamma$ instead of $4\gamma$ above. We choose such form of Riesz kernel simply to match the spacial regularity of the solution with the one from the bounded domain case. 

We recall that $G_{\alpha}(t,x)=\displaystyle\int_{\mathbb{R}}e^{ix\xi-t|\xi|^{\alpha}}\dif\xi$ is the fundamental solution of 
$\partial_tG_{\alpha}(t,x) = -\left(-\Delta\right)^{\alpha/2}G_{\alpha}(t,x)$. Consequently, the mild solution to \eqref{eq:fracHeatR} is defined as 
\begin{align}\label{eq:mildSol}
	X_t(x)=\sigma\int_0^t\int_{\mathbb{R}} G_{\alpha}\left(\theta(t-s),x-z\right)W^{\gamma}\left(\dif s,\dif z\right),
\end{align}
where the above integral is a Wiener integral with respect to the Gaussian noise $W^\gamma$. For details, see for instance \cite[Section~3]{LototskyRozovsky2017Book} and \cite[Section~2]{Dalang1999}. 

In the context of statistical inference, SPDEs similar to \eqref{eq:fracHeatR} were recently considered in \cite{KhalilTudor2019} and \cite{KhalilTudor2019b}.

\begin{proposition}
	For $m\in \mathbb{N}\cup \{0\}$, we have that 
	$$
	\sup_{t\in[0,T],x\in\mathbb{R}}\mathbb{E}\left|\partial_x^m X_t(x)\right|^2<\infty\quad \mbox{for every}\quad T>0,
	$$
	if and only if $1+2m<\alpha+4\gamma$. 
	In particular,  $\partial_x^m X_t(x)$ is well-defined for $x\in\R$, $0\leq t\leq T$ and  $1+2m<\alpha+4\gamma$. 
	
\end{proposition}
\begin{proof}
	Without loss of generality, we fix $\theta=\sigma=1$ for simplicity. Then, for each $0\leq t\leq T$ and $x\in \mathbb{R}$, we have
	\begin{align*}
		&\mathbb{E}\left|\partial_x^m X_t(x)\right|^2
		= (2\pi)^{-1}\int_0^t \int_{\mathbb{R}} |\xi|^{2m-4\gamma}e^{-2s|\xi|^{\alpha}}\dif \xi \dif s\\
		&=\int_0^t s^{(4\gamma-2m-1)/\alpha}\dif s\int_{\mathbb{R}}e^{-2|\xi|^{\alpha}}|\xi|^{2m-4\gamma}\dif \xi<\infty,
	\end{align*}
	if and only if $(4\gamma-2m-1)/\alpha>-1$, which is equivalent to $1+2m<\alpha+4\gamma$.
\end{proof}

Define the following remainder term:
\begin{align}
	R(x) &= \sigma\int_t^\infty\int_\R \partial_x^m\left(G_\alpha\left(\theta (t-s), z\right) - G_\alpha\left(\theta (t-s), x-z\right)\right) W^\gamma(\dif s, \dif z).
\end{align}
Note that the remainder term decays exponentially fast in Fourier space and is therefore smooth in space for each $t>0$.

The next result is based on \cite[Proposition 4.6]{KhalilTudor2019b}. 

\begin{proposition}\label{prop:fullSpaceReg}
	For $t>0$, the process $\partial_x^mX_t$ has the same distribution as a perturbed fBM of the form $c_{\alpha,\gamma,m}\displaystyle\frac{\sigma}{\sqrt{\theta}}B^{2\gamma+\alpha/2-1/2 - m} + R$, provided that $2\gamma+\alpha/2-1/2 - m\in (0, 1)$, where $c_{\alpha,\gamma,m}^2:=(2\pi)^{-1}\int_{\mathbb{R}}\left(1-\cos(\xi)\right)|\xi|^{2m-4\gamma-\alpha}\dif \xi$ and $R\in C^{\infty}(\mathbb{R})$ almost surely.
\end{proposition}
\begin{proof}
	For every $x\in \mathbb{R}$ and a fixed $t>0$, we set $	v(x):=\partial_x^m X_t(x)-R(x)$. Then, for $x,y\in \mathbb{R}$,
	\begin{align*}
		\mathbb{E}\left|v(x)-v(y)\right|^2\nonumber
		&=\frac{\sigma^2}{\pi}\int_0^{\infty}\int_{\mathbb{R}}|\xi|^{2m-4\gamma}\left(1-\cos(\xi(x-y))\right) e^{-2\theta s|\xi|^{\alpha}}\dif \xi \dif s\\
		&= \frac{\sigma^2}{2\pi\theta}(x-y)^{\alpha+4\gamma-1-2m}\int_{\mathbb{R}}\left(1-\cos(\xi)\right)|\xi|^{2m-4\gamma-\alpha}\dif \xi.
	\end{align*}
	We note that, by the assumption $0<(\alpha+4\gamma-1)/2 - m<1$,
	\begin{align*}
		\int_{\mathbb{R}}\left(1-\cos(\xi)\right)|\xi|^{2m-4\gamma-\alpha}\dif \xi&\lesssim \int_{|\xi|>1} |\xi|^{2m-4\gamma-\alpha}d\xi + \int_{|\xi|\leq 1} |\xi|^{2m-4\gamma-\alpha+2}d\xi<\infty.
	\end{align*}
	This implies that $v$ is a fractional Brownian motion with Hurst index $\frac{\alpha+4\gamma-1}{2}-m$. The smoothness property  $R\in C^{\infty}(\mathbb{R})$ follows from \cite[Proposition 4.6]{KhalilTudor2019b}. This concludes the proof. 
\end{proof}

The following result on estimation of drift $\theta$ or volatility $\sigma$ of fractional heat equation \eqref{eq:fracHeatR} follows immediately from Theorem~\ref{thm:itfBmCLT} in conjunction with Proposition~\ref{prop:fullSpaceReg}. 

\begin{theorem}
	Let $m\in \mathbb{N}_0$ and $0<H<1$ such that $m+H=s^*=2\gamma+\alpha/2-1/2$. Let $M,q\in \mathbb{N}$, and assume that either $M=m+1$ with $H<1/2$ or $M\geq m+2$. Then, we have, as $N\to \infty$,
	\begin{align}
	\sqrt{N}\left(V_{q,M,s^*,N}\left(X_t\right)-c_{\alpha,\gamma,m}^q\tau_q\mu_{M,s^*}^{q/2}\left(\frac{\sigma}{\sqrt{\theta}}\right)^q\right)\overset{d}\to\mathcal{N}\left(0,c_{\alpha,\gamma,m}^{2q}\sigma_{q,M,s^*}^2\mu_{M,s^*}^q\left(\frac{\sigma}{\sqrt{\theta}}\right)^{2q}\right).
	\end{align}
	Moreover,
	\begin{itemize}
		\item[(i)] If $\theta$ is known, then $\widetilde{\sigma}_N^{q,M}:=c_{\alpha,\gamma,m}^{-1}\tau_q^{-1/q}\mu_{M,s^*}^{-1/2}\sqrt{\theta}V_{q,M,s^*,N}(X_t)^{1/q}$ is an asymptotically normal estimator for $\sigma$;
		\item[(ii)] If $\sigma$ is known, $\widetilde{\theta}_N^{q,M}:=c_{\alpha,\gamma,m}^{2}\tau_q^{2/q}\mu_{M,s^*}\sigma^2V_{q,M,s^*,N}(X_t)^{-2/q}$ is an asymptotically normal estimator for $\theta$.
	\end{itemize}
\end{theorem}

We conclude this section with several clarifying remarks on the class of considered SPDEs in this section. The choice of Riesz kernel was primarily prompted by \cite{KhalilTudor2019b} that considers same equations. 
In particular this allows to have a direct compassion of the results obtained in this paper and those from \cite{KhalilTudor2019,KhalilTudor2019b}. A careful reader will also notice that working with Riesz kernel, which is characterized by its Fourier transform $\mathcal{F}K_{\gamma}(\xi)=|\xi|^{-4\gamma}$,  is technically convenient. On the other hand, such correlation structure of the noise limits $\gamma \in(0,1/4)$, thus limiting the range of regularity of the solution in spatial component (as described above). To overcome this, but also to be on par with SPDEs from Section~\ref{sec:SPDE-BDom}, one can replace the Riesz kernel with Bessel kernel  with Fourier transform $\mathcal{F}K^B_\gamma(\xi) = (1+|\xi|^2)^{-2\gamma}$, for any $\gamma>0$. This case indeed can be addressed, and results similar to those from Section~\ref{sec:SPDE-BDom} can be obtained. For the sake of brevity, we shortly sketch the main arguments of the proof. For simplicity, let us also assume 
that the drift operator $-\theta(-\Delta)^{\alpha/2}$ is substituted by $-\theta(I-\Delta)^{\alpha/2}$ in \eqref{eq:fracHeatR}. First, we note that for $0<\gamma<1/4$, there exists a positive definite kernel $K^{R/B}_\gamma$ such that $\mathcal{F}K^{R/B}_\gamma = \mathcal{F}K_\gamma-\mathcal{F}K^B_\gamma$. Let $m\in\R$ such that
$\gamma':=2\gamma+\alpha/2-1/2-m\in(0,1)$, let $\widetilde W^{\gamma'}$ a centered Gaussian field with covariance kernel $K^{R/B}_{\gamma'}$, independent of $W^\gamma$. Then, similarly to \cite[Proposition 4.6]{KhalilTudor2019b}, one can prove that the increments of $(I-\Delta)^{m/2}u-R^{(1)}-R^{(2)}$ are the increments of a fractional Brownian motion, where
\begin{align*}
	R^{(1)}(x) &= \sigma\int_t^\infty\int_\R(I-\Delta)^{m/2}(G_\alpha(\theta(t-s),z)-G_\alpha(\theta(t-s), x-z))W^\gamma(\mathrm{d}z, \mathrm{d}s), \\
	R^{(2)}(x) &= \sigma\int_0^\infty\int_\R(G_\alpha(\theta(t-s), z)-G_\alpha(\theta(t-s), x-z))\widetilde W^{\gamma''}(\mathrm{d}z, \mathrm{d}s).
\end{align*}
Then, for $m\in 2\N$,  we have that $(I-\Delta)^{m/2}X$ is a linear combination of $\partial_{x}^{2m'}X$, where $0\leq m'\leq m/2$, so $V_{q,M,s,N}(X)=V_{q,M,s,N}(J^m\partial_x^mu)=V_{q,M,s,N}(J^m(I-\Delta)^{m/2}X)$. Furthermore, $J^m(I-\Delta)^{m/2}X$ behaves like $J^mB^H$ with $H=\gamma'$, up to a perturbation by $R^{(1)}$ and $R^{(2)}$. Consequently, similar statements concerning consistency and rate of convergence of the $\Delta$--power variation can be made. However, the central limit theorem does not transfer since $R^{(2)}$ is not arbitrarily smooth. Similar to the bounded domain, the asymptotic normality property of the corresponding estimators for $\theta$ and $\sigma$ remains an open question.

The emphasis that the extension of the results from linear to nonlinear equations of the form \eqref{eq:intro-main-bdd} via a splitting argument depends on spatial regularity properties of the solution $X$ to \eqref{eq:fracHeatR}. In contrast to the case of bounded domains, the covariance operator as given by the Riesz (or Bessel) kernel is not of trace class, so $X$ will not belong to $L^2(\R)$ or any higher-order Sobolev space derived from $L^2(\R)$. Instead, we believe that suitably chosen weighted Sobolev spaces can help to mitigate this issue. To the best of our knowledge, this has not been investigated systematically in the literature.

\section{Numerical example}\label{sec:examples}

In this section we illustrate the theoretical results of Section~\ref{sec:SPDE-BDom} via numerical simulations, by considering  the stochastic heat equation
\begin{align}
	\mathrm{d}X_t = \theta\Delta X_t\mathrm{d}t + \sigma(-\Delta)^{-\gamma}\mathrm{d}W_t,
\end{align}
with initial condition $X_0=0$ on $\mathcal{D}=[0,1]$ with Dirichlet boundary conditions. We take the true  values of the parameters $\theta, \sigma=1$. 
As far as the smoothing parameter $\gamma$, we consider the following representative cases $\gamma\in\{0.0, 0.375, 0.5, 0.625\}$, which correspond to the regularity level $s^*=2\gamma+1/2\in\{0.5, 1.25, 1.5, 1.75\}$. To numerically simulate a path, we use the Fourier series decomposition of the solution \eqref{eq:SPDE-Fourier} by taking $N_0=1\times 10^4$ eigenmodes, and each eigenmode is numerically simulated by the Euler implicit scheme with  temporal stepsize $\delta t=1\times 10^{-8}$. Correspondingly, the solution is computed at $N_0+1$ uniformly spaced spatial grid points with step size $h= 1\times10^{-4}$. 

\begin{figure}
	\includegraphics[width=0.5\textwidth]{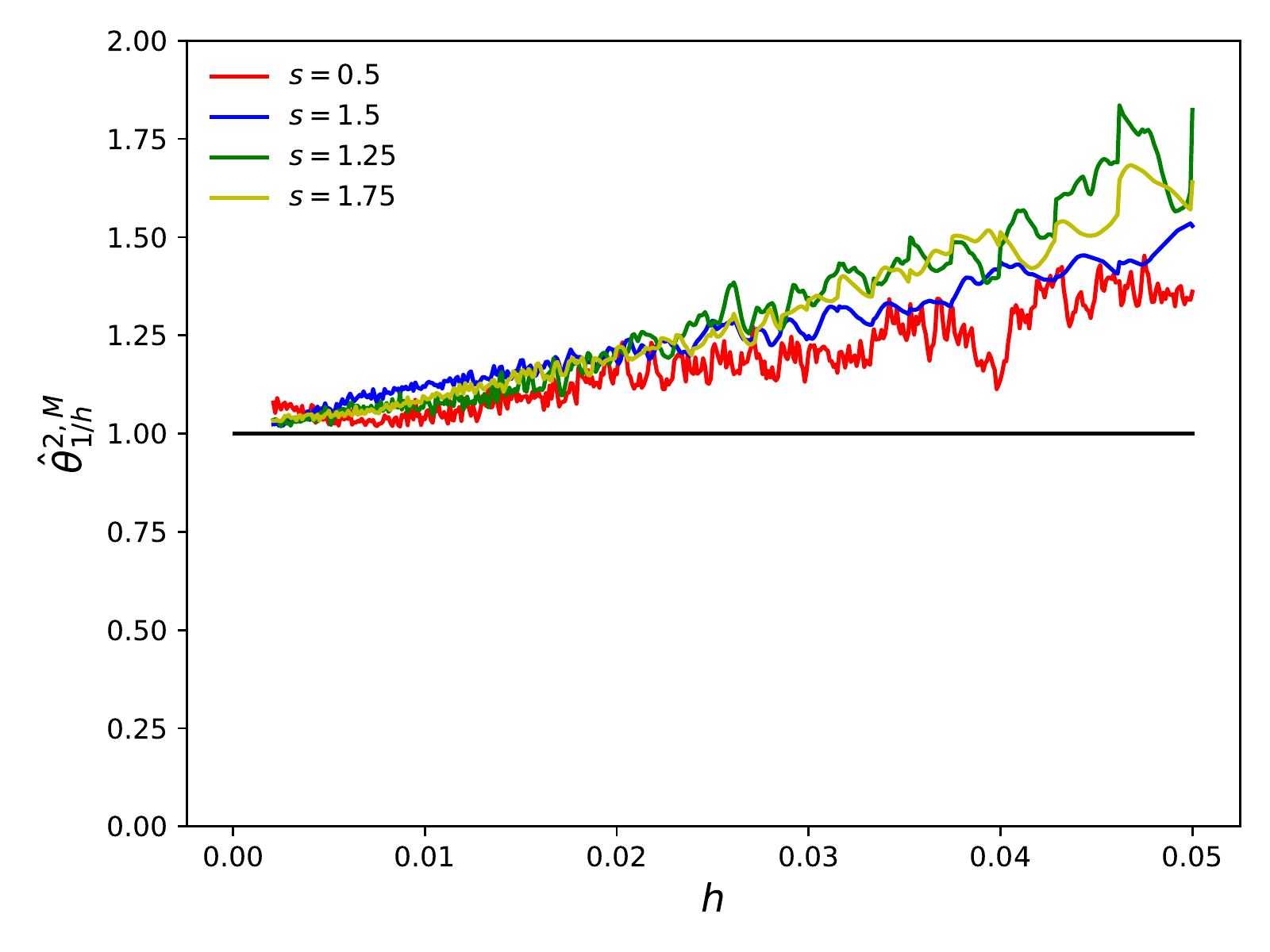}
	\includegraphics[width=0.5\textwidth]{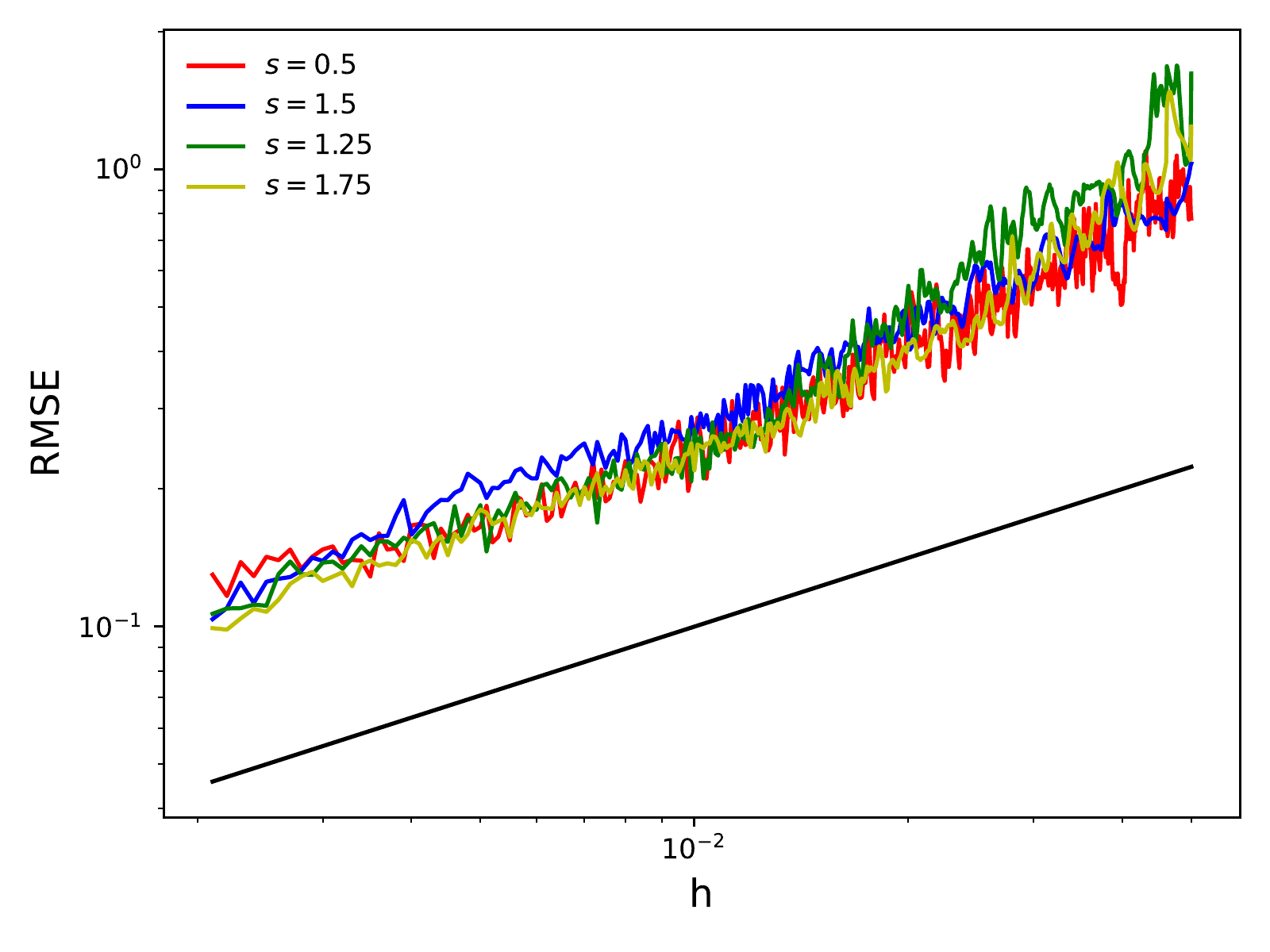}
	\includegraphics[width=0.5\textwidth]{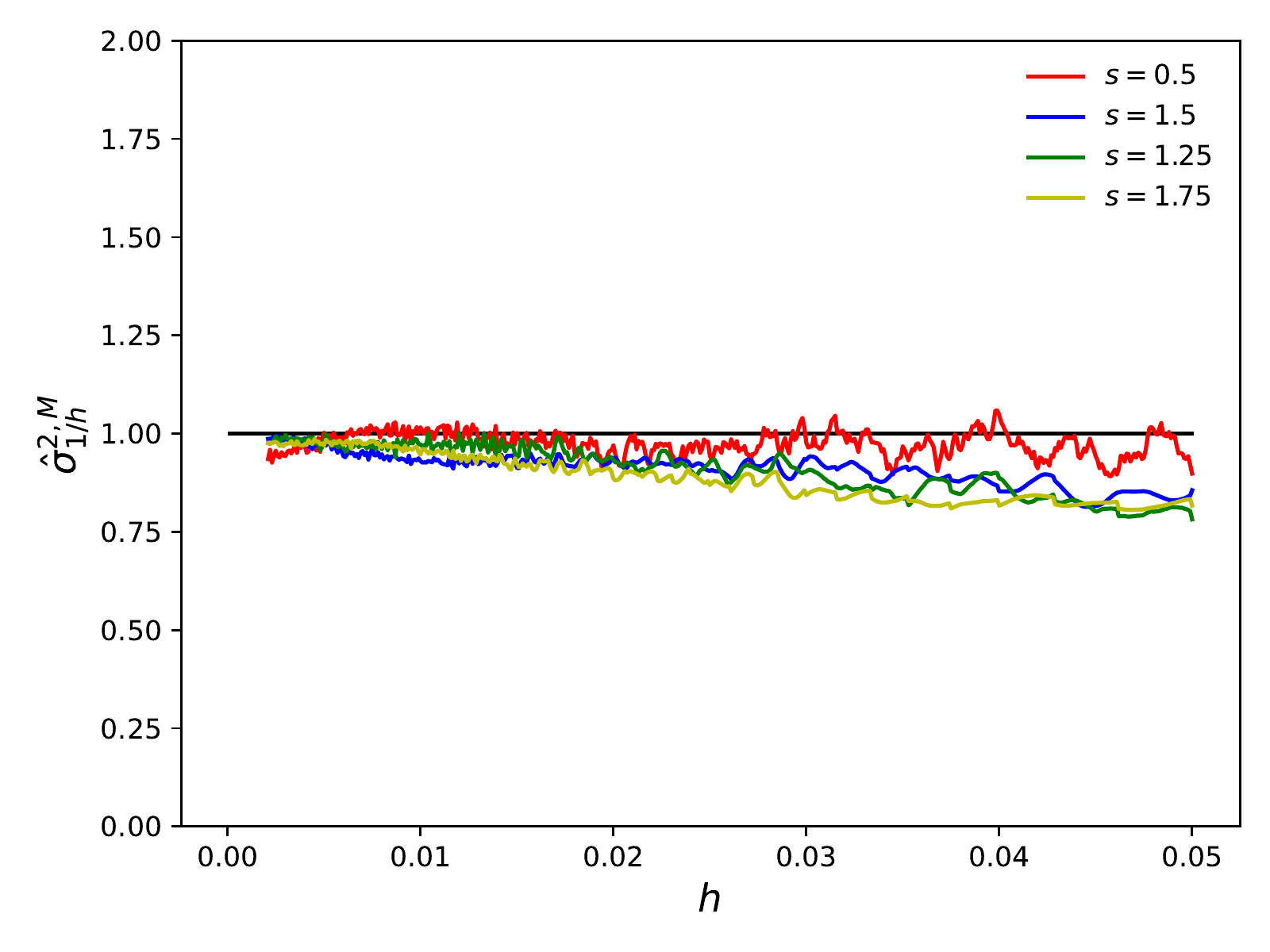}
	\includegraphics[width=0.5\textwidth]{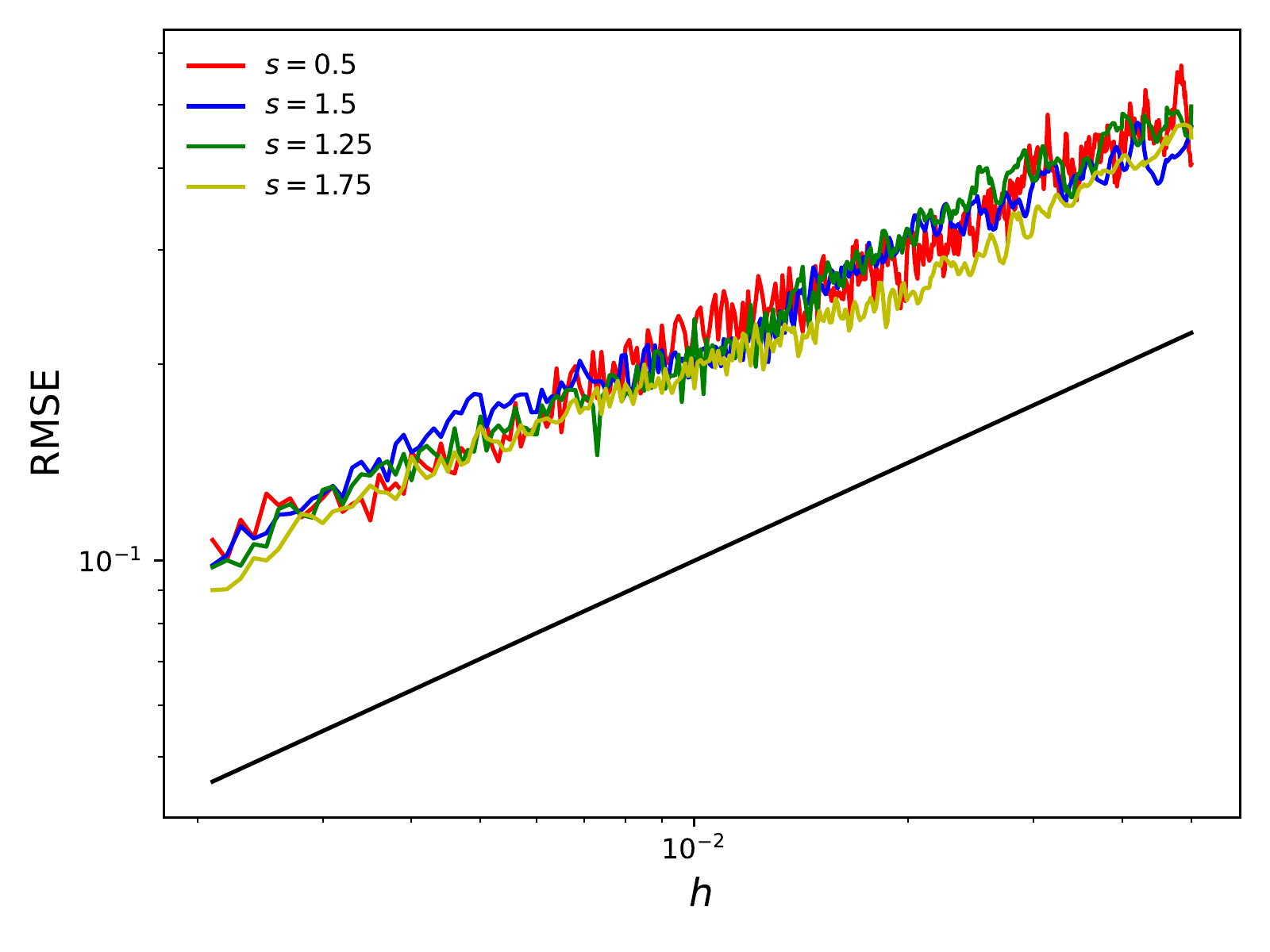}
	\caption{Estimation of  $\theta$ (top row) and $\sigma^2$ (bottom row). Left panel: the average of $100$ Monte Carlo estimates as function of spatial sampling resolution $h$. The solid black line corresponds to the true value $1.0$. Right panel:  The RMSE (root mean square error) as function of $h$. The black line corresponds to the theoretical convergence rate $h^{1/2}$. 
		\label{fig:numerics}
	}
\end{figure}

Next, we assume that the solution $X$ is observed at time $T=1$ on spatial grid points belonging to the interval $[a,b]$, with $a=0.2$, $b=0.8$. We apply Theorem~\ref{th:SPDE-parmEst}, with 
$q=2$ and $M=\lceil s^*\rceil + 2$, to estimate  one of the parameters $\mu$ or $\sigma^2$ assuming that the second one is known. For each set of the parameters, we perform these evaluations on 100 Monte Carlo sample paths of the solution. The average values of the estimates as function of step size $h$ are displayed in Figure~\ref{fig:numerics}, left panel. Clearly, the estimators converge to the true value (horizontal solid line $\theta=1$ and $\sigma^2=1$), as the mesh size gets smaller. Moreover, as shown in Figure~\ref{fig:numerics}, right panel, the root mean square error of the estimators behaves as $h^{1/2}$, confirming the theoretical rate of converges of the proposed estimators, regardless of the order of regularity $s^*$ of the solution. Similar results were obtained for various sets of the parameters. Finally, while not shown here, we remark that the results from Section~\ref{sec:fBM} were also confirmed via numerical simulations. 

The numerical computations were performed using programing language Python. The source
code is available from the authors upon request.

\section*{Acknowledgments}
IC acknowledges partial support from the National Science Foundation grant DMS-1907568.
The research of GP has been funded by Deutsche Forschungsgemeinschaft (DFG) - SFB1294/1 - 318763901. 
GP thanks the Illinois Institute of Technology for the hospitality during a research visit, where this project has been initiated.

\appendix

\section{Appendix} 

{\small  
For reader's convenience, we recall a useful asymptotic result of Hermite polynomials of a stationary Gaussian sequence.

\begin{theorem}\cite[Theorem 7.2.4 Breuer-Major Theorem]{NourdinPeccati2012}\label{thm:Breuer-Major}\label{th:Breuer-Major}
	Let $Y=\{Y_k\}_{k\in \mathbb{Z}}$ be a centered stationary Gaussian sequence with unit variance, and $f(x)=\displaystyle\sum_{q=d}^{\infty} a_qH_q(x),\ a_q\in \mathbb{R},$
	where $H_q$ is the $q$-th Hermite polynomial.
	Assume that 
	\begin{equation}\label{eq:sumEll}
		\sum_{\ell\in \mathbb{Z}}\left|\rho(\ell)\right|^d<\infty, 
	\end{equation}
	where $\rho(\ell)=\mathbb{E}\left(Y_0Y_{\ell}\right),\ \ell\in \mathbb{Z}$.
	Then,
	$$
	\lim_{N\to \infty}\frac{1}{\sqrt{N}}\sum_{k=1}^N f(Y_k) \overset{d}=\mathcal{N}\left(0,\sum_{q=d}^{\infty} q!a_q^2\sum_{\ell\in \mathbb{Z}} \rho(\ell)^q\right).
	$$
\end{theorem}
}

\newcommand{\etalchar}[1]{$^{#1}$}
\def\cprime{$'$}


\end{document}